\input amstex
\documentstyle{amsppt}
\input epsf
\magnification=\magstep1 
\pageheight{9.0truein}
\pagewidth{6.5truein}
\parindent=20pt
\NoBlackBoxes

\def\sD{\Cal D}

\NoBlackBoxes
\topmatter

\author
A.N. Dranishnikov
\endauthor

\title
On generalized amenability
\endtitle
\abstract 
There is a word metric $d$ on countably generated free group $\Gamma$
such that $(\Gamma,d)$ does not admit a coarse uniform imbedding into
a Hilbert space.
\endabstract


\address Pennsylvania State University, Department of Mathematics, 
218 McAllister Building, University Park, PA, 16802, USA
\endaddress

\subjclass Primary 20H15 
\endsubjclass

\email  dranish\@math.psu.edu
\endemail

\keywords amenable group, coarsely uniform embedding 
\endkeywords

\endtopmatter

\document
\head \S1 Introduction \endhead

A discrete countable group $G$ is called {\it amenable} if there exists a
left invariant {\it mean} on $G$, i.e. a positive finitely additive, finite
measure $\mu$. Clearly, that $\mu(g)=0$ for
all $g\in G$. Equivalently, a group $G$ is amenable if its natural action
on the Stone-\v{C}ech compactification $\beta G$ admits an invariant
measure. In [Gr] M. Gromov introduced the notion of  an {\it a}-\rm{T}-
{\it menable} group as a group $G$ which admits a proper isometric action
on the Hilbert space $l_2$.

The Novikov higher signature conjecture was known for some classes of amenable groups for
many years. Recently Higson and Kasparov [H-K] proved it for 
all amenable groups and for a-T-menable
groups. Then G. Yu [Y] proved it for more general
class of groups, we call it {\it Y-amenable} groups. A group $G$ is called Y-amenable if it admits a coarsely uniform
embedding as defined in [Gr] into a Hilbert space.

In the case of genuine amenability there is the Folner Criterion [Fo],[Gr]
which allows to establish amenability of a group in terms of the growth 
function of an exhausting family of compact sets in a group.
In [Y] Yu introduced his Property A
(we do not define it here), which serves as a distant analog of Folner
property. After analyzing the Property A Higson and Roe [H-R] introduced
a new notion of amenability.

\proclaim{Definition} 
A discrete countable group $G$ is called  Higson-Roe amenable if its action
on the Stone-\v{C}ech compactification $\beta G$ is topologically amenable.
\endproclaim
An action of $G$ on a compact space $X$ is {\it topologically amenable}
[A-D-R] if there is a sequence of continuous maps $b^n:X\to P(G)$
to the space of probability measures on $G$ such that for every $g\in G$,
$\lim_{n\to\infty}\sup_{x\in X}\|gb^n_x-b_{gx}^n\|_1=0$.
Here a measure $b^n_x=b^n(x)$ is considered as a function $b^n_x:G\to[0,1]$ and
$\|\ \|_1$ is the $l_1$-norm.

\proclaim{Assertion 1}
A discrete countable group $G$ is Higson-Roe amenable if it admits a
topologically amenable action on some compact metrizable space $X$.
\endproclaim
\demo{Proof}
The proof in one direction is given in Proposition 2.3 of [H-R]. The other 
implication follows from
countability of $G$ and the Schepin Spectral Theorem [Ch].\qed
\enddemo
We note that the trivial 
action of the classiscal amenable groups on a one-point space is
topologically amenable. Also all hyperbolic groups are acting 
on their Gromov boundaries topologically amenable [Ad], [A-D-R]. Still
there is no example of a countable group which is not Higson-Roe 
amenable. In this paper we present an example of countable
group which is not Y-amenable.  

\head \S2 Coarsely uniform embeddings \endhead

A map $f:X\to Y$ between metric spaces is called a 
{\it coarsely uniform embedding} if there are functions
$\rho_1,\rho_2;[0,\infty)\to[0,\infty)$ with 
$\lim_{t\to\infty}\rho_i(t)=\infty$ such that

$\rho_1(d_X(x,x'))\le d_Y(f(x),f(x'))\le\rho_2(d_X(x,x'))$ for all $x,x'\in X$.

{\bf 1. Higson-Roe-Yu Embedding Theorem.} The following theorem is
due to Higson-Roe and Yu [H-R],[Y].
\proclaim{Theorem 1}
 A finitely generated Higson-Roe amenable group $G$ admits a coarsely 
uniform embedding into the Hilbert space for a word metric on $G$.
\endproclaim
Every set of generators $S$ of a group $G$ defines a word metric
$d_S$ on $G$. If a group is finitely generated we assume that $S$ is
finite. Any two such metrics generated by two finite
sets are quazi-isometric. The following fact is well-known.
\proclaim{Assertion 2}
Let $\Gamma$ be a finitely generated subgroup of a finitely
generated group $G$, then the inclusion is a coarsely uniform
embedding.
\endproclaim
\demo{Proof}
Let $S$ be a set of generators of $\Gamma$ and let $T$ be a
set of generators of $G$. Without loss of generality we may
assume that $S\subset T$. Then $d_T(x,y)=\|x^{-1}y\|_T\le
\|x^{-1}y\|_S=d_S(x,y)$. Thus, $\rho_2(t)=t$. We define
$\rho_1(t)=\min\{\|w\|_T\mid w\in\Gamma,\|w\|_S\ge t\}$.
Assume that $\rho_1$ is bounded. Then there are a constant $R$ and 
a sequence $w_i\in\Gamma$ with $\|w_i\|_S\ge i$ and $\|w_i\|_T\le R$.
This contradicts with the fact that a $R$-ball in $G$ is finite.
Note that $\rho_1(d_S(x,y))\le\|x^{-1}y\|_T=d_T(x,y)$.\qed
\enddemo

{\bf 2. Modified Enflo's Metric Spaces.}
We define metric spaces $M_n$ which are adaptations for asymptotic geometry
of Enflo's spaces [En]. Let $N_n=\{0,1,2,\dots, 2^{n+1}-1\}$ with
metric $|x-y| \mod 2^{n+1}$. We define
$M_n=(N_n)^{2n^n}$ as the product of $2n^n$ copies of $N_n$ with the
metric $d(a,b)=\max_i\{|a_i-b_i|\mod 2^{n+1}\}$ where $a=\{a_i\}$ and
$b=\{b_i\}$.

A pair of points $(a,b)$ in $M_n$ is called an $m$-segment if the coordinates of $a$ and $b$
are different in exactly $2n^{n-m}$ positions and $|a_i-b_i|=2^m$ 
if $a_i\ne b_i$.
\proclaim{Proposition 1}
For any two $m$-segments $(a,b)$ and $(a',b')$ in $M_n$ there is an isometry
$h:M_n\to M_n$ with $h(a)=a'$ and $h(b)=b'$ such that $h$ takes $k$-segments to
$k$-segments for any $k$.
\endproclaim
\demo{Proof}
First we consider a permutation $\sigma:\{1,2,\dots,2n^n\}\to\{1,2,\dots,2n^n\}$
which establishes a bijection between coordinate spaces for which $a_i=b_i$ and
$a'_i=b'_i$. Then for every $i$ we consider an isometry
$h_i:N_n\to N_n$ taking $(a_{\sigma(i)},b_{\sigma(i)})$ to $(a'_i,b'_i)$.
Such an isometry exists, since either $|a_{\sigma(i)}-b_{\sigma(i)}|=2^m=|a'_i-b'_i|$
or $|a_{\sigma(i)}-b_{\sigma(i)}|=0=|a'_i-b'_i|$.
The family $\{h_i\}$ defines an isometry $\bar h:M_n\to M_n$. We define
$h=\bar h\circ\bar\sigma$ where $\bar\sigma:M_n\to M_n$ is defined by the formula
$\bar\sigma(x_1,\dots, x_{2n^n})=(x_{\sigma(1)},\dots,x_{\sigma(2n^n)})$.
Then $(h(a))_i=(\bar h\circ\bar\sigma(a))_i=h_i(\bar\sigma(a))_i)=
h_i(a_{\sigma(i)})=a_i'$. Thus, $h(a)=a'$. Similarly, one can check that
$h(b)=b'$.\qed
\enddemo

Following Enflo [En], by a {\it double $n$-simplex} in a space $M$ 
we call a set $2n+2$ points
$D_n=\{u_1,\dots,u_{n+1},v_1,\dots, v_{n+1}\}$ ,\  $u_i,v_j\in M$.
Pairs $(u_i,u_j)$ and $(v_i,v_j)$ are called {\it edges} of $D_n$ and pairs 
$(u-i,v_j)$ are called {\it connecting lines}.
\proclaim{Proposition 2}
For any $m$, $1\le m< n$, there exists a double $(n-1)$-simplex
$D^m_{n-1}\subset M_n$ such that all edges are $m$-segments and 
all connecting lines are $(m-1)$-segments.
\endproclaim
\demo{Proof}
Let $I_m=\{1,\dots,n^{n-m+1}\}$. Let $J^m_1,\dots,J^M_n$ be the partition 
of $I_M$ in $n$ equal parts $J^m_1=\{1,\dots,n^{n-m}\}$,
$J^m+2=\{n^{n-m}+1,\dots,2n^{n-m}\}$, ..., $ J^M_n=\{(n-1)n^{n-m}+1,
\dots,n^{n-m+1}\}$.
We define $u_k,v_k\in M_n$ as follows:

$$(u_k)_i=\cases
2^m,&if\  i\in J_k^m;\\
2^{m-1},&if\ i\in I_m+n^{n-m+1};\\
0,& otherwise.
\endcases
$$
and

$$
(v_k)_i=\cases
2^m & if \ i\in J^m_k+n^{n-m+1};\\
2^{m-1} &if \ i\in I_m;\\
0 & otherwise.
\endcases
$$

Since $u_k$ and $u_l$ for $k\ne l$ differ at $2n^{n-m}$ positions and
$|(u_k)_i-(u_l)_i|=2^m$ at those positions,
 all $u$-edges in the corresponding double $(n-1)$-simplex are $m$-segments.
 Similarly, all $v$-edges are $m$-segments. Since $u_k$ and $v_l$ are distinct
in $2n^{n-m+1}$ coordinates with $|(u_k)_i-(v_l)_i|=2^{m-1}$, every 
 connecting line $(u_k,v_l)$ is an $(m-1)$-segment.\qed
 \enddemo

{\bf 3. Obstruction to Embedding}
The following proposition is well-known, it can be extracted from [En].
\proclaim{Proposition 3} For every double $n$-simplex in the Hilbert space
the inequality $\Sigma c_{\alpha}^2\ge\Sigma s^2_{\beta}$ holds where 
$c_{\alpha}$ runs through the length of connecting lines and $s_{\beta}$
runs through the length of edges.
\endproclaim
\demo{Proof}
First we proof this inequality for a double simplex in the real line.
The equality

$\Sigma_{1\le k,l\le n}(u_k-v_l)^2-\Sigma_{1\le k<l\le n}(u_k-u_l)^2-
\Sigma_{1\le k<l\le n}(v_k-v_l)^2$

$=(\Sigma_{1\le k\le n}u_k-\Sigma_{1\le l\le n}b_l)^2$

implies the inequality

$\Sigma_{1\le k,l\le n}(u_k-v_l)^2\ge\Sigma_{1\le k<l\le n}(u_k-u_l)^2+
\Sigma_{1\le k<l\le n}(v_k-v_l)^2$ 

which is exactly the inequality
$\Sigma c_{\alpha}^2\ge\Sigma s^2_{\beta}$.

Since $\|u_k-u_l\|^2=\Sigma_i((u_k)_i-(v_l)_i)^2$, we obtain the required 
inequality by adding up corresponding inequalities for $i$-th coordinates.
\qed
\enddemo
\proclaim{Theorem 2}
Assume that a metric space $X$ contains an isometric copies of $M_n$ 
for all $n$ Then $X$ cannot be coarsely uniformly embedded in the
Hilbert space.
\endproclaim
\demo{Proof}
Assume the contrary. Let $f:X\to l_2$ be a coarsely uniform embedding.
let $\rho_1$ and $\rho_2$ be corresponding functions. Since $\rho_1\to\infty$,
there is $m$ such that $\rho_1(2^m)>2\sqrt{e}\rho_2(1)$.
We consider a double $(n-1)$-simplex $D^m_{n-1}\subset M_n\subset X$ for any
$n>m$. Denote by $\bar f((a,b))=\|f(a)-f(b)\|$. Then Proposition 3 implies the
inequality 

\

$\Sigma_{c\in C(D^m_{n-1})}\bar f(c)^2\ge
\Sigma_{s\in E(D^m_{n-1})}\bar f(s)$. 

\

Here $C(D^m_{n-1})$ denotes the set of all connecting
lines and $E(D^m_{n-1})$ denotes the set of all edges. 

Let $\sD$ be the set of all double
$(n-1)$ simplices in $M_n$ isomorphic to $D^m_{n-1}$. Then

\

$\Sigma_{c\in C(D),D\in\sD}\bar f(c)^2\ge\Sigma_{s\in E(D),D\in\sD}\bar f(s)^2$.

\

This inequality can be written as

\

$\Sigma_{c\in C(D),D\in\sD}\frac{\bar f(c)^2}{n^2card(\sD)}\ge
\Sigma_{s\in E(D),D\in\sD}\frac{\bar f(s)^2}{n^2card(\sD)}=\frac{n-1}{n}
\Sigma_{s\in E(D),D\in\sD}\frac{\bar f(s)^2}{n(n-1)card(\sD)}$.

\

Let $S_m$ denote the set of all $m$-segments in $M_n$. By the Proposition 1
all $m$-segments in $M_n$ are equal. It means that every $m$-segment $c$
is a connecting line for the same number of double simplices from $\sD$ and
every $m-1$-segment is an edge of the same number of double simplices from
$\sD$.
Since the number of connecting edges in a double $(n-1)$-simplex is $n$,
the expression $\Sigma_{c\in C(D),D\in\sD}\frac{\bar f(c)^2}{n^2card(\sD)}$ is
the arithmetic mean.
Because of symmetry
the arithmetic mean 
$\Sigma_{c\in C(D),D\in\sD}\frac{\bar f(c)^2}{n^2card(\sD)}$
can be computed as $\Sigma_{c\in S_{m-1}}\frac{\bar f(c)^2}{card(S_{m-1})}$.
By a similar reason the arithmetic mean
$\Sigma_{s\in E(D),D\in\sD}\frac{\bar f(s)^2}{n(n-1)card(\sD)}$ can be computed
as $\Sigma_{c\in S_{m}}\frac{\bar f(c)^2}{card(S_{m})}$.
Thus, we have an inequality

\

$(1+\frac{1}{n-1})\bar g_{m-1}\ge \bar g_m$ where $\bar g_k=
\Sigma_{c\in S_k}\frac{\bar f(c)^2}{card(S_k)}$.

\

Iterate this inequality to obtain the following

$(1+\frac{1}{n-1})^{n-1}\bar g_o\ge \bar g_{n-1}$. Hence $e\bar g_0\ge\bar g_{n-1}$.
Then 

$\sqrt{e}\rho_2(1)\ge\sqrt{e}\sup_{c\in S_0}\bar 
f(c)\ge\inf_{c\in S_{n-1}}\bar f(c)\ge
\rho_1(2^{n-1})\ge\rho_1(2^m)>2\sqrt{e}\rho_2(1)$.

The contradiction completes the proof.\qed
\enddemo 

{\bf 4. A group which is not Y-amenable.}

For every modified Enflo's space $M_n$ we consider the graph $G_n$
whose vertices are points of $M_n$ and two vertices $a$ and $b$ are joint
by an edge if and only if $d(a,b)=1$ in $M_n$. 
Note that $G_n$ is connected.
Let $G$ be an infinite wedge of all $G_n$. We define a path metric
on $G$ such that any two vertices joined by an edge are on distance one.
We define a countable infinitely generated group $\Gamma$ as follows.
Fix an orientation on all edges of $G$. Then
all edges of $G$ are the generators of $\Gamma$ and all loops are the relations.
\proclaim{Theorem 3}
The group $\Gamma$ is not Y-amenable.
\endproclaim
\demo{Proof}
Fix a metric on $\Gamma$ defined by the above set of generators, then
$G$ is isometrically imbedded in $\Gamma$. By Theorem 2 $\Gamma$ does not admit
a coarsely uniform embedding into the Hilbert space. \qed
\enddemo
It can be shown that the group $\Gamma$ is infact a 
 free group generated by edges of a maximal tree in $G$.

I am very grateful to M. Gromov, N. Higson, J. Roe and M. Sapir for
valuable discussions and remarks.

\enddocument